\documentclass[11pt]{article}
\usepackage{graphicx,psfrag}
\usepackage{amssymb,amsfonts,amsthm}
\usepackage{amsmath,amsthm,amssymb}
\usepackage{amsfonts}

\newcommand{\calr}{{\mathcal R}}
\newtheorem{theorem}{Theorem}[section]
\newtheorem{lemma}[theorem]{Lemma}

\newtheorem{proposition}[theorem]{Proposition}
\theoremstyle{definition}
\newtheorem{definition}[theorem]{Definition}
\newtheorem{definitions}[theorem]{Definitions}

\newtheorem{remark}[theorem]{Remark}
\newtheorem{prob}[theorem]{Problem}

\newcommand{\call}{{\mathcal L}}
\newcommand{\calg}{\mathcal{G}}
\newcommand{\dist}{{\mathrm{dist}}}
\newcommand{\ddd}{{\mathcal D}}
\newcommand{\uuu}{{\mathcal U}}
\newcommand{\eee}{{\mathcal E}}
\newcommand{\nn}{{\mathcal N}}

\newcommand {\Z}{\mathbb{Z}}            
\newcommand {\R}{\mathbb{R}} 
\newcommand {\C}{\mathbb{C}} 
\newcommand {\q}{\mathfrak q} 

\newcommand {\me}{\medskip}

\newcommand {\iv}{^{-1}}

\begin{document}

\title{Relatively Hyperbolic Groups with Rapid Decay Property}
\author{Cornelia Dru\c{t}u and Mark Sapir\thanks{The first author
is grateful to the CNRS of France for granting her the
d\'el\'egation CNRS status in the 2003-2004 academic year. The
research of the second author was supported in part by the NSF grant
DMS 0245600.}}
\date{}
\maketitle

\me

\begin{abstract}
We prove that a finitely generated group $G$ that is (strongly)
relatively hyperbolic with respect to a collection of finitely
generated subgroups $\{ H_1,\dots , H_m\}$ has the Rapid Decay
property if and only if each $H_i\, ,\, i=1,2,\dots ,m$, has the
Rapid Decay property.
\end{abstract}

\section{Introduction}

Throughout the whole paper, unless stated otherwise, $G$ denotes a finitely generated group, and $1$ denotes the neutral element in $G$.


The $L^2$-norm of a function $x$ in $l^2(G)$ is denoted by $\| x
\|$.

A group $G$ satisfies \textit{the Rapid Decay property} (RD
property, in short) if the space of rapidly decreasing
functions on $G$ with respect to some length-function is inside the
reduced $C^*$-algebra of $G$ (see Section \ref{prel} for the precise
definition). The RD property is relevant to the Novikov conjecture
\cite{CM} and to the Baum-Connes conjecture \cite{L}.

The main result of the paper is the following.

\begin{theorem}[Theorem \ref{main}, Proposition \ref{centre}]\label{imain}
Let $G$ be a group which is (strongly) relatively hyperbolic with respect
 to some finitely generated subgroups
$\{ H_1,H_2,\dots , H_m\}$. Then $G$ has the RD property if and only
if $\{ H_1,\dots , H_m\}$ have the RD property.
\label{m}\end{theorem}

The ``only if" part in Theorem \ref{imain} follows from the more
general statement that a subgroup of a group that has the property
RD also has RD with respect to the induced length-function
\cite[Proposition 2.1.1]{J}. The ``if" part is more difficult. A
proof of it is given in Section \ref{sproof}. In fact we prove that
the statement of the theorem still holds if we replace (strong)
relative hyperbolicity by a natural weaker version of it, that we
call (*)-relative hyperbolicity (see Definition \ref{star}). At the
end of the paper, we shall show how to further weaken property (*)
to get an even wider class of groups with RD (see Remark
\ref{twostars}).

The argument in the proof of Theorem \ref{imain} is similar in
spirit to the arguments used in \cite{RRS}, \cite{L1} and \cite{T}.
The idea, in all the mentioned arguments, is to reduce the proof of
the inequality (\ref{ineq0}) to the case when the convolution is
performed on some ``easier" geodesic triangles. In our case, the
``easier" triangles are simply the triangles contained in left
cosets $gH_i$; we can pass from general triangles to triangles in
$gH_i$ because of the (*)-relative hyperbolicity.

Some particular cases of Theorem \ref{imain} have been proven
before:

\begin{itemize}

\item All hyperbolic groups satisfy the RD property \cite{J,H,C} (this is a
particular case of Theorem \ref{imain} with $m=1$, $H_1=\{1\}$).

\item The amalgamated product of two groups $A$ and $B$ with finite
amalgamated subgroup $F$ satisfies RD provided $A$ and $B$ satisfy
RD \cite{J} (take $G=A*_F B, H_1=A, H_2=B$).

\item A group $G$ relatively hyperbolic with respect to the subgroups $H_1,\dots ,H_m$ ,
 has the RD property provided that $H_1,\dots ,H_m$ have polynomial growth
\cite{CR}. The fact that polynomial growth implies RD follows from
the definition of RD \cite{J}.
\end{itemize}

Not all groups have RD: a group that contains an amenable subgroup
of superpolynomial growth (with respect to the word metric of the
whole group) does not have RD \cite{J}. Note that it is the only
known obstruction to RD.

\section{Preliminaries}
\subsection{Property RD}\label{prel}

Recall that a length-function on a group $G$ is a map $L$ from $G$
to the set of non-negative real numbers $\R_+$ satisfying:
\begin{itemize}
    \item[(1)] $L(gh)\leq L(g)+ L(h)$  for all $g,h\in G$ ;
   \item[(2)] $L(g)=L(g\iv)$  for all $g\in G$ ;
    \item[(3)] $L(1)=0$.
\end{itemize}

\me

\noindent \textit{Notations}: We denote by
$\overline{B}_L (r)$ the $L$-ball of radius $r$, that is the set $\{
g\in G \mid L (g)\leq r\}$ and by $S_L (r)$ the $L$-sphere of radius
$r$, that is the set $\{ g\in G \mid L (g)= r\}$.

\me

We say that a length-function $L_1\colon G\to \R_+$
\textit{dominates} another length-function $L_2\colon G\to \R_+$ if
there exist $a,b\in \R_+$ such that $L_2\leq aL_1+b$. If $L_1$
dominates $L_2$ and $L_2$ dominates $L_1$ then $L_1$ and $L_2$ are
said to be \textit{equivalent}.

\begin{remark}\label{balls}
If $L_1$ dominates $L_2$ then there exists $c\in \R_+$ such that
for every $r\geq 1$ we have $\overline{B}_{L_1} (r)\subset
\overline{B}_{L_2} (cr)$.
\end{remark}

If $G$ is a finitely generated group, the length-functions
corresponding to two finite generating sets are equivalent. All such
length-functions are called {\textit{word length-functions.}}

\begin{lemma}[\cite{J}, Lemma 1.1.4]\label{word}
If $G$ is a finitely generated group then any word length-function
dominates all length-functions on $G$.
\end{lemma}

\me

We first give the traditional analytic version of the RD property with respect to a length-function $L$, and then present a more geometric one.  For every $s\in \R$ {\it the
Sobolev space of order $s$ with respect to $L$} is the set $H_L^s
(G)$ of functions $\phi$ on $G$ such that the function $(1+L)^s\phi$
is in $l^2 (G)$. The space of {\it{rapidly decreasing functions on
$G$ with respect to $L$}} is the set $H_L^\infty (G)=\bigcap_{s\in
\R }H_L^s (G)$.

The {\it{group algebra of $G$}}, denoted by $\C G$, is the set of
functions with finite support on $G$, i.e. it is the set of formal
linear combinations of elements of $G$ with complex coefficients. We
denote by $\R_+ G$ its subset consisting of functions taking values
in $\R_+$.

With every element $g\in G$ we can associate the linear {\em
convolution operator} $\phi\mapsto g*\phi$ on $l^2(G)$, where
$$g*\phi(h)=\phi(g\iv h).$$ This is just the left regular representation of $G$ on $l^2(G)$, it can be extended to a representation of $\mathbb{C}G$ on
$l^2(G)$ by linearity. This representation is faithful and every
convolution operator induced by an element of $\C G$ is bounded.
Therefore we can identify $\C G$ with a subspace in the space of
bounded operators $\mathbf{B}(l^2(G))$ on $l^2(G)$. For every $x\in
\C G$ we denote by $\|x\|_*$ its operator norm, that is
$$
\| x\|_* = \sup \{ \|x*\phi\| \; ;\; \|\phi\| =1\}\, .
$$

The closure $C_r^* (G)$ of $\C G$ in the operator norm is called
{\it{the reduced $C^*$-algebra of $G$}}.


\begin{definition}
The group $G$ is said to have {\it{the RD property with respect to the length-function $L$}} if the inclusion of $\C G$ into $C_r^*(G)$ extends to a continuous inclusion of $H_L^\infty
(G)$ into $C_r^*(G)$.
\end{definition}

We recall an equivalent way of defining RD. The
following result is a slight modification of Proposition 1.4 in
\cite{CR}.

\begin{lemma}\label{equiv} Let $G$ be a discrete group and let $L$ be a length-function on it. The following statements are equivalent:
\begin{itemize}
    \item[(i)] The group $G$ has the RD property with respect to
    $L$.
    \item[(ii)] There exists a polynomial $P$ such that for every
    $r>0$, every $x\in \R_+G$ such that $x$ vanishes outside
    $\overline{B}_L(r)$, and every $\phi\in l^2(G)$ such that $\phi(G)\subseteq
    \R_+$, we have
\begin{equation}\label{ineq0}
\| x * \phi\| \leq P(r)\, \| x \| \cdot \|\phi \|\, .
\end{equation}
\end{itemize}
\end{lemma}

\proof In \cite[Proposition 1.4]{CR} it is proved that (i) is
equivalent to (ii) for $x\in \R_+ G$ and $\phi\in l^2(G)$. We prove
that if (ii) is satisfied for $\phi$ with $\phi(G)\subseteq \R_+$
then it is satisfied for every $\phi\in l^2(G)$. Let $\phi\in
l^2(G)$. Then we can write
$$
\phi=\phi_1 - \phi_2 +i (\phi_3 - \phi_4),
$$ where $\phi_i$ take values in $\R_+$ and $\|\phi\|^2=\sum_{i=1}^4 \|\phi_i\|^2$.
We have the inequalities
$$
\| x*\phi\| \leq \sum_{i=1}^4 \| x*\phi_i\| \leq P(r)
\|x\|\sum_{i=1}^4 \|\phi_i\|\leq 2 P(r) \|x\| \cdot \|\phi\|.
$$
\endproof

\begin{definition}
A group is said to \textit{have the RD property} if it satisfies the
RD property with respect to some length-function.
\end{definition}

\begin{remark}
Lemma \ref{word} and Remark \ref{balls} imply that if a finitely
generated group satisfies RD with respect to some length-function,
then it satisfies RD with respect to a word length-function (for
some generating set). Hence, a finitely
generated group satisfies RD if and only if it
satisfies RD with respect to every word length-function.
\end{remark}

In the case of finitely generated groups, it suffices to show (\ref{ineq0}) for functions with finite support, as shown by the following lemma.

\me

\noindent \textit{Notation}: For a function $f\in l^2(G)$ and a constant $p \geq 0$, $f_p$
denotes the function which coincides with $f$ on $S_L(p)$ and which
vanishes outside $S_L(p)$.

\me

\begin{lemma}\label{equiv1} Let $G$ be a finitely generated group and let $L$ be a word length-function on it. The following statements are equivalent:
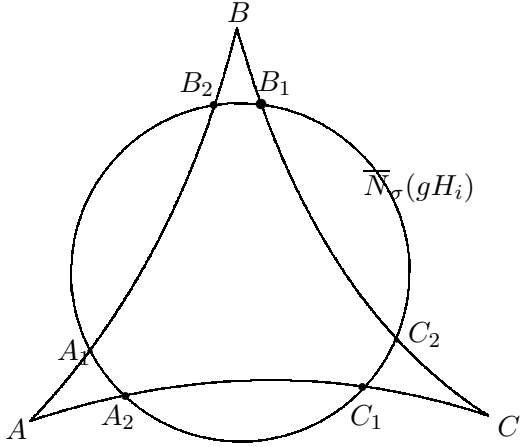
\begin{itemize}
    \item[(i)] The group $G$ has the RD property.
\item[(ii)] There exists a polynomial $P$
such that for every
    $r_1,r_2\ge 0$, every $p\in [| r_1-r_2 | \, ,\, r_1+r_2]$
    every $x\in \R_+G$
    with support in
    $S_L(r_1)$, and every $y\in \R_+ G$
    with support in $S_L(r_2)$,
\begin{equation}\label{ineq}
    \| (x * y)_p \| \leq P(r)\, \| x \| \cdot \|y \|\, .
\end{equation}
\end{itemize}
\end{lemma}

\proof The equivalence is true according to the argument in the proof of Theorem  5, \cite[$\S III.5.\alpha $]{C}.\endproof

\me

For details on the RD property we refer to \cite{C}, \cite{J},
\cite{L1}, \cite{CR}.

\subsection{(*)-relative hyperbolicity}

\noindent \textit{Notation}: For a subset $Y$ in a metric space we
denote by $\overline{\nn}_\delta (Y)$ the closed tubular
neighborhood of $Y$, that is $\{ x\mid \dist (x,Y) \leq \delta
\}$.

\me

\begin{definition}\label{star}

Let $G$ be a group and let $H_1,...,H_m$ be subgroups in $G$. We say
that $G$ is (*)-{\textit{relatively hyperbolic with respect to}}
$H_1,...,H_m$ if there exists a finite generating set $S$ of $G$,
and two constants $\sigma$ and $\delta$ such that the following
property holds:

\begin{itemize}
\item[(*)] For every geodesic triangle $ABC$ in the Cayley graph of $G$ with
respect to $S$, there exists a coset $gH_i$ such that
$\overline{\nn}_\sigma(gH_i)$ intersects each of the sides of the
triangle, and the entrance (resp. exit) points $A_1, B_1, C_1$
(resp. $B_2, C_2, A_2$) of the sides $[A_,B], [B,C], [C,A]$ in
$\overline{\nn}_\sigma(gH_i)$ satisfy $$\dist(A_1,A_2)<\delta,\,
\dist(B_1, B_2)<\delta,\, \dist(C_1, C_2)<\delta\, .$$
\end{itemize}
\end{definition}

Note that $H_1,...,H_m$ need not be finitely generated.

\begin{figure}
\centering
\unitlength .7mm 
\linethickness{0.4pt}
\ifx\plotpoint\undefined\newsavebox{\plotpoint}\fi 
\begin{picture}(105.5,87)(0,0)

\qbezier(14.5,9.25)(42.38,38.38)(53.75,84)
\qbezier(53.75,84)(68.38,31.38)(101.5,10.25)
\qbezier(101.5,10.25)(59.5,24.25)(14.5,9.25)
\put(86.65,37.5){\line(0,1){1.26}}
\put(86.63,38.76){\line(0,1){1.258}}
\multiput(86.55,40.02)(-.0308,.3135){4}{\line(0,1){.3135}}
\multiput(86.43,41.27)(-.02872,.20804){6}{\line(0,1){.20804}}
\multiput(86.26,42.52)(-.03159,.17722){7}{\line(0,1){.17722}}
\multiput(86.03,43.76)(-.03369,.15386){8}{\line(0,1){.15386}}
\multiput(85.77,44.99)(-.03176,.12194){10}{\line(0,1){.12194}}
\multiput(85.45,46.21)(-.03319,.10964){11}{\line(0,1){.10964}}
\multiput(85.08,47.42)(-.0317,.0916){13}{\line(0,1){.0916}}
\multiput(84.67,48.61)(-.032746,.083838){14}{\line(0,1){.083838}}
\multiput(84.21,49.78)(-.033606,.076991){15}{\line(0,1){.076991}}
\multiput(83.71,50.94)(-.032291,.066719){17}{\line(0,1){.066719}}
\multiput(83.16,52.07)(-.032943,.061769){18}{\line(0,1){.061769}}
\multiput(82.57,53.18)(-.033479,.05725){19}{\line(0,1){.05725}}
\multiput(81.93,54.27)(-.032297,.050571){21}{\line(0,1){.050571}}
\multiput(81.25,55.33)(-.032697,.047027){22}{\line(0,1){.047027}}
\multiput(80.53,56.37)(-.033014,.043722){23}{\line(0,1){.043722}}
\multiput(79.77,57.37)(-.033256,.040629){24}{\line(0,1){.040629}}
\multiput(78.97,58.35)(-.033429,.037722){25}{\line(0,1){.037722}}
\multiput(78.14,59.29)(-.03354,.034984){26}{\line(0,1){.034984}}
\multiput(77.27,60.2)(-.034885,.033643){26}{\line(-1,0){.034885}}
\multiput(76.36,61.07)(-.037624,.03354){25}{\line(-1,0){.037624}}
\multiput(75.42,61.91)(-.040531,.033375){24}{\line(-1,0){.040531}}
\multiput(74.45,62.71)(-.043625,.033142){23}{\line(-1,0){.043625}}
\multiput(73.44,63.48)(-.046931,.032835){22}{\line(-1,0){.046931}}
\multiput(72.41,64.2)(-.050476,.032445){21}{\line(-1,0){.050476}}
\multiput(71.35,64.88)(-.057151,.033647){19}{\line(-1,0){.057151}}
\multiput(70.26,65.52)(-.061672,.033125){18}{\line(-1,0){.061672}}
\multiput(69.15,66.12)(-.066624,.032487){17}{\line(-1,0){.066624}}
\multiput(68.02,66.67)(-.072086,.031718){16}{\line(-1,0){.072086}}
\multiput(66.87,67.18)(-.083741,.032992){14}{\line(-1,0){.083741}}
\multiput(65.7,67.64)(-.091506,.031969){13}{\line(-1,0){.091506}}
\multiput(64.51,68.05)(-.10954,.03351){11}{\line(-1,0){.10954}}
\multiput(63.3,68.42)(-.12185,.03212){10}{\line(-1,0){.12185}}
\multiput(62.08,68.74)(-.13668,.03035){9}{\line(-1,0){.13668}}
\multiput(60.85,69.02)(-.17713,.03211){7}{\line(-1,0){.17713}}
\multiput(59.61,69.24)(-.20796,.02933){6}{\line(-1,0){.20796}}
\multiput(58.37,69.42)(-.3134,.0317){4}{\line(-1,0){.3134}}
\put(57.11,69.54){\line(-1,0){1.258}}
\put(55.85,69.62){\line(-1,0){2.52}}
\put(53.33,69.63){\line(-1,0){1.258}}
\put(52.08,69.56){\line(-1,0){1.254}}
\multiput(50.82,69.44)(-.24975,-.03373){5}{\line(-1,0){.24975}}
\multiput(49.57,69.27)(-.17731,-.03107){7}{\line(-1,0){.17731}}
\multiput(48.33,69.05)(-.15396,-.03324){8}{\line(-1,0){.15396}}
\multiput(47.1,68.79)(-.12203,-.0314){10}{\line(-1,0){.12203}}
\multiput(45.88,68.47)(-.10974,-.03287){11}{\line(-1,0){.10974}}
\multiput(44.67,68.11)(-.091692,-.03143){13}{\line(-1,0){.091692}}
\multiput(43.48,67.7)(-.083934,-.032499){14}{\line(-1,0){.083934}}
\multiput(42.31,67.25)(-.077089,-.033379){15}{\line(-1,0){.077089}}
\multiput(41.15,66.75)(-.066814,-.032095){17}{\line(-1,0){.066814}}
\multiput(40.01,66.2)(-.061866,-.032761){18}{\line(-1,0){.061866}}
\multiput(38.9,65.61)(-.057348,-.03331){19}{\line(-1,0){.057348}}
\multiput(37.81,64.98)(-.050666,-.032148){21}{\line(-1,0){.050666}}
\multiput(36.75,64.3)(-.047123,-.032558){22}{\line(-1,0){.047123}}
\multiput(35.71,63.59)(-.043819,-.032885){23}{\line(-1,0){.043819}}
\multiput(34.7,62.83)(-.040726,-.033136){24}{\line(-1,0){.040726}}
\multiput(33.72,62.04)(-.037821,-.033318){25}{\line(-1,0){.037821}}
\multiput(32.78,61.2)(-.035083,-.033437){26}{\line(-1,0){.035083}}
\multiput(31.87,60.33)(-.032496,-.033498){27}{\line(0,-1){.033498}}
\multiput(30.99,59.43)(-.033651,-.037525){25}{\line(0,-1){.037525}}
\multiput(30.15,58.49)(-.033494,-.040432){24}{\line(0,-1){.040432}}
\multiput(29.34,57.52)(-.03327,-.043527){23}{\line(0,-1){.043527}}
\multiput(28.58,56.52)(-.032972,-.046834){22}{\line(0,-1){.046834}}
\multiput(27.85,55.49)(-.032593,-.05038){21}{\line(0,-1){.05038}}
\multiput(27.17,54.43)(-.032124,-.0542){20}{\line(0,-1){.0542}}
\multiput(26.53,53.35)(-.033306,-.061574){18}{\line(0,-1){.061574}}
\multiput(25.93,52.24)(-.032683,-.066528){17}{\line(0,-1){.066528}}
\multiput(25.37,51.11)(-.031929,-.071992){16}{\line(0,-1){.071992}}
\multiput(24.86,49.96)(-.033238,-.083644){14}{\line(0,-1){.083644}}
\multiput(24.4,48.78)(-.032238,-.091412){13}{\line(0,-1){.091412}}
\multiput(23.98,47.6)(-.03102,-.10032){12}{\line(0,-1){.10032}}
\multiput(23.6,46.39)(-.03247,-.12175){10}{\line(0,-1){.12175}}
\multiput(23.28,45.18)(-.03075,-.13659){9}{\line(0,-1){.13659}}
\multiput(23,43.95)(-.03263,-.17703){7}{\line(0,-1){.17703}}
\multiput(22.77,42.71)(-.02994,-.20787){6}{\line(0,-1){.20787}}
\multiput(22.59,41.46)(-.0327,-.3133){4}{\line(0,-1){.3133}}
\put(22.46,40.21){\line(0,-1){1.257}}
\put(22.38,38.95){\line(0,-1){1.26}}
\put(22.35,37.69){\line(0,-1){1.26}}
\put(22.37,36.43){\line(0,-1){1.258}}
\put(22.43,35.17){\line(0,-1){1.255}}
\multiput(22.55,33.92)(.033,-.24985){5}{\line(0,-1){.24985}}
\multiput(22.72,32.67)(.03054,-.1774){7}{\line(0,-1){.1774}}
\multiput(22.93,31.42)(.03279,-.15406){8}{\line(0,-1){.15406}}
\multiput(23.19,30.19)(.03104,-.12213){10}{\line(0,-1){.12213}}
\multiput(23.5,28.97)(.03255,-.10983){11}{\line(0,-1){.10983}}
\multiput(23.86,27.76)(.031161,-.091784){13}{\line(0,-1){.091784}}
\multiput(24.26,26.57)(.032252,-.084029){14}{\line(0,-1){.084029}}
\multiput(24.72,25.39)(.033152,-.077187){15}{\line(0,-1){.077187}}
\multiput(25.21,24.24)(.031898,-.066908){17}{\line(0,-1){.066908}}
\multiput(25.76,23.1)(.032579,-.061962){18}{\line(0,-1){.061962}}
\multiput(26.34,21.98)(.033141,-.057446){19}{\line(0,-1){.057446}}
\multiput(26.97,20.89)(.033599,-.053298){20}{\line(0,-1){.053298}}
\multiput(27.64,19.83)(.032419,-.047219){22}{\line(0,-1){.047219}}
\multiput(28.36,18.79)(.032756,-.043916){23}{\line(0,-1){.043916}}
\multiput(29.11,17.78)(.033016,-.040823){24}{\line(0,-1){.040823}}
\multiput(29.9,16.8)(.033207,-.037918){25}{\line(0,-1){.037918}}
\multiput(30.73,15.85)(.033334,-.035181){26}{\line(0,-1){.035181}}
\multiput(31.6,14.93)(.033402,-.032594){27}{\line(1,0){.033402}}
\multiput(32.5,14.05)(.035987,-.032462){26}{\line(1,0){.035987}}
\multiput(33.44,13.21)(.040334,-.033613){24}{\line(1,0){.040334}}
\multiput(34.41,12.4)(.043429,-.033398){23}{\line(1,0){.043429}}
\multiput(35.4,11.64)(.046737,-.03311){22}{\line(1,0){.046737}}
\multiput(36.43,10.91)(.050284,-.032741){21}{\line(1,0){.050284}}
\multiput(37.49,10.22)(.054105,-.032283){20}{\line(1,0){.054105}}
\multiput(38.57,9.57)(.061476,-.033487){18}{\line(1,0){.061476}}
\multiput(39.68,8.97)(.066432,-.032879){17}{\line(1,0){.066432}}
\multiput(40.81,8.41)(.071898,-.032141){16}{\line(1,0){.071898}}
\multiput(41.96,7.9)(.083546,-.033484){14}{\line(1,0){.083546}}
\multiput(43.13,7.43)(.091316,-.032507){13}{\line(1,0){.091316}}
\multiput(44.31,7.01)(.10023,-.03131){12}{\line(1,0){.10023}}
\multiput(45.52,6.63)(.12166,-.03283){10}{\line(1,0){.12166}}
\multiput(46.73,6.3)(.1365,-.03116){9}{\line(1,0){.1365}}
\multiput(47.96,6.02)(.17693,-.03315){7}{\line(1,0){.17693}}
\multiput(49.2,5.79)(.20778,-.03056){6}{\line(1,0){.20778}}
\multiput(50.45,5.61)(.3132,-.0336){4}{\line(1,0){.3132}}
\put(51.7,5.47){\line(1,0){1.257}}
\put(52.96,5.39){\line(1,0){1.26}}
\put(54.22,5.35){\line(1,0){1.26}}
\put(55.48,5.36){\line(1,0){1.259}}
\put(56.73,5.43){\line(1,0){1.255}}
\multiput(57.99,5.54)(.24994,.03226){5}{\line(1,0){.24994}}
\multiput(59.24,5.7)(.17749,.03002){7}{\line(1,0){.17749}}
\multiput(60.48,5.91)(.15416,.03234){8}{\line(1,0){.15416}}
\multiput(61.72,6.17)(.12222,.03068){10}{\line(1,0){.12222}}
\multiput(62.94,6.48)(.10993,.03222){11}{\line(1,0){.10993}}
\multiput(64.15,6.83)(.09953,.03346){12}{\line(1,0){.09953}}
\multiput(65.34,7.23)(.084123,.032005){14}{\line(1,0){.084123}}
\multiput(66.52,7.68)(.077284,.032925){15}{\line(1,0){.077284}}
\multiput(67.68,8.17)(.071189,.033683){16}{\line(1,0){.071189}}
\multiput(68.82,8.71)(.062057,.032397){18}{\line(1,0){.062057}}
\multiput(69.93,9.3)(.057543,.032972){19}{\line(1,0){.057543}}
\multiput(71.03,9.92)(.053397,.033442){20}{\line(1,0){.053397}}
\multiput(72.1,10.59)(.047314,.03228){22}{\line(1,0){.047314}}
\multiput(73.14,11.3)(.044012,.032627){23}{\line(1,0){.044012}}
\multiput(74.15,12.05)(.04092,.032896){24}{\line(1,0){.04092}}
\multiput(75.13,12.84)(.038016,.033095){25}{\line(1,0){.038016}}
\multiput(76.08,13.67)(.035279,.03323){26}{\line(1,0){.035279}}
\multiput(77,14.53)(.032692,.033306){27}{\line(0,1){.033306}}
\multiput(77.88,15.43)(.032568,.035891){26}{\line(0,1){.035891}}
\multiput(78.73,16.37)(.033731,.040235){24}{\line(0,1){.040235}}
\multiput(79.54,17.33)(.033526,.043331){23}{\line(0,1){.043331}}
\multiput(80.31,18.33)(.033247,.046639){22}{\line(0,1){.046639}}
\multiput(81.04,19.35)(.032889,.050188){21}{\line(0,1){.050188}}
\multiput(81.73,20.41)(.032442,.05401){20}{\line(0,1){.05401}}
\multiput(82.38,21.49)(.033667,.061377){18}{\line(0,1){.061377}}
\multiput(82.99,22.59)(.033074,.066335){17}{\line(0,1){.066335}}
\multiput(83.55,23.72)(.032352,.071803){16}{\line(0,1){.071803}}
\multiput(84.07,24.87)(.03373,.083447){14}{\line(0,1){.083447}}
\multiput(84.54,26.04)(.032775,.09122){13}{\line(0,1){.09122}}
\multiput(84.96,27.22)(.03161,.10014){12}{\line(0,1){.10014}}
\multiput(85.34,28.43)(.03319,.12156){10}{\line(0,1){.12156}}
\multiput(85.67,29.64)(.03156,.13641){9}{\line(0,1){.13641}}
\multiput(85.96,30.87)(.03367,.17684){7}{\line(0,1){.17684}}
\multiput(86.19,32.11)(.03117,.20769){6}{\line(0,1){.20769}}
\multiput(86.38,33.35)(.0276,.2505){5}{\line(0,1){.2505}}
\put(86.52,34.61){\line(0,1){1.257}}
\put(86.61,35.86){\line(0,1){1.638}}
\put(12,7.75){\makebox(0,0)[cc]{$A$}}
\put(54.25,87){\makebox(0,0)[cc]{$B$}}
\put(105.5,8.25){\makebox(0,0)[cc]{$C$}}
\put(23,22.25){\makebox(0,0)[cc]{$A_1$}}
\put(61,73.5){\makebox(0,0)[cc]{$B_1$}}
\put(78.5,10){\makebox(0,0)[cc]{$C_1$}}
\put(46.25,73.25){\makebox(0,0)[cc]{$B_2$}}
\put(89.5,25.75){\makebox(0,0)[cc]{$C_2$}}
\put(31.25,10){\makebox(0,0)[cc]{$A_2$}}
\put(49.5,69.25){\circle*{1.41}} \put(58.5,69.5){\circle*{1.8}}
\put(84.25,24.75){\circle*{1.12}} \put(77.75,15.75){\circle*{1.41}}
\put(32.75,14){\circle*{1.5}} \put(26,22.5){\circle*{1}}
\put(88.5,54){\makebox(0,0)[cc]{$\overline{N}_\sigma(gH_i)$}}
\end{picture}

\end{figure}

The following statement is proved in \cite{DS}.

\begin{proposition}[\cite{DS}, Corollary 8.14, Lemma 8.19]\label{centre}
Every group $G$ that is (strongly) relatively hyperbolic with
respect to a family of finitely generated subgroups $H_1,...,H_m$ is
(*)-relatively hyperbolic with respect to these subgroups.
\end{proposition}

Note that the converse statement is not true. For example, every
hyperbolic group $G$ is obviously (*)-relatively hyperbolic with
respect to any subgroup $H<G$. But it is well known that for strong
relative hyperbolicity to hold $H$ must be quasiconvex.

\begin{prob} Is it true that every group that is (*)-relatively
hyperbolic with respect to certain subgroups $H_1,...,H_m$ is also
strongly relatively hyperbolic with respect to some subgroups
$H_1',...,H_n'$ such that each $H_i'$ is inside some $H_j$?
\end{prob}

Note that there exists a version of the definition of the strong
relative hyperbolicity which also makes sense without $H_1,...,H_m$
being finitely generated \cite{Osin}.

\begin{prob}
Can one remove the hypothesis that the subgroups $H_1,...,H_m$ are
finitely generated from the Proposition \ref{centre} ?
\end{prob}

\section{Proof of the main result}\label{sproof}

\begin{theorem}\label{main}
Let $G$ be a group which is (*)-relatively hyperbolic with respect
to the subgroups $\{ H_1,H_2,\dots , H_m\}$. Then $G$ has the RD
property if and only if $\{ H_1,\dots , H_m\}$ have the RD property
with respect to the length-function induced by a word
length-function on $G$.
\end{theorem}

\proof Only the ``if" part needs a proof. Consider a finite generating set $S$
of $G$ with respect to which the property (*) is satisfied, and consider the
word length-function $L$ on $G$ corresponding to $S$. By hypothesis and the fact
that two word length-functions are equivalent, we may suppose that
$H_1,...,H_m$ have RD with respect to $L$.


\me

\noindent {\textit{Convention}}: In what follows we fix two
arbitrary positive numbers $r_1,r_2$ and we fix $p \in
[|r_1-r_2|\, ,\, r_1+r_2]$.

\me

Let $P_i(r)$ be the polynomial given by Lemma \ref{equiv}, (ii),
for the group $H_i,\, i=1,2,\dots ,m$, and the length function
$L$, and let $P(r)=1+\sum_{i=1}^m P_i(r+2\kappa )^2$, where
$\kappa$ is a constant to be chosen later. We shall prove that for
every $x$ with support in $S_L(r_1)$ and every $y$ with support in
$S_L(r_2)$, we have the inequality
\begin{equation}\label{2}
\|(x*y)_p\|^2 \leq Q(r_1)P(r_1) \|x \|^2 \cdot \|y\|^2 \, ,
\end{equation} where
$Q$ is a polynomial of degree $3$. That will imply inequality
(\ref{ineq}).

For every $g\in S_L (r)$ we choose one geodesic $\q_g$ joining it to
$1$. We thus obtain a set of geodesics of length $r$ indexed by
elements in $S_L(r)$. We denote this set of geodesics by $\calg(r)$ and
we identify it with the set $S_L(r)$.

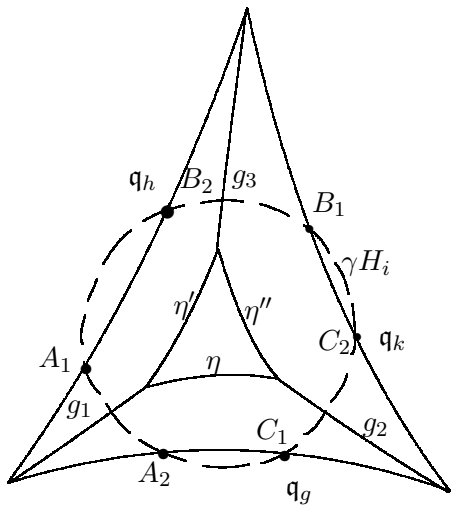
\begin{figure}[!ht]
\centering
\unitlength .7mm 
\linethickness{0.4pt}
\ifx\plotpoint\undefined\newsavebox{\plotpoint}\fi 
\begin{picture}(91.25,93.75)(0,20)
\qbezier(7.25,23.5)(35,61.38)(52.75,113.75)
\qbezier(52.75,113.75)(65,56.13)(91.25,22)
\qbezier(91.25,22)(50.25,36.75)(7.25,23.5)
\put(33,81){\makebox(0,0)[cc]{$\q_h$}}
\put(80.5,50.5){\makebox(0,0)[cc]{$\q_k$}}
\put(62.75,21.75){\makebox(0,0)[cc]{$\q_g$}}
\qbezier(7.5,23.5)(20,31.88)(33.5,41.75)
\qbezier(33.5,41.75)(41.13,52.13)(47.25,68)
\qbezier(47.25,68)(52.5,49.88)(58.75,43.25)
\qbezier(58.75,43.25)(47.75,45.5)(33.75,41.75)
\qbezier(58.25,43.5)(74.38,31.88)(91,21.75)
\qbezier(47,67.5)(48.25,80.88)(52.5,112.75)
\put(52.5,81.25){\makebox(0,0)[cc]{$g_3$}}
\put(77.25,34){\makebox(0,0)[cc]{$g_2$}}
\put(46.5,45.75){\makebox(0,0)[cc]{$\eta$}}
\put(41,57.25){\makebox(0,0)[cc]{$\eta'$}}
\put(55,56.25){\makebox(0,0)[cc]{$\eta''$}}
\put(21,37.5){\makebox(0,0)[cc]{$g_1$}}
\qbezier(25.5,65.25)(28.38,70.13)(30.75,71.5)
\qbezier(33,73.25)(36.13,75.25)(40.75,76.25)
\qbezier(44,76.75)(48.5,77.63)(52,77)
\qbezier(55.25,76.5)(58.5,75.63)(60.75,74.25)
\qbezier(63.5,72.5)(67.13,70.13)(68.25,67.25)
\qbezier(69.5,65.25)(71.25,63.13)(72,59.5)
\qbezier(72.75,57.5)(73.88,48.25)(72.5,48)
\qbezier(71.75,44.25)(71.25,41.88)(68.75,39)
\qbezier(67.75,37.25)(65.38,33)(63.5,32.75)
\qbezier(62.25,30.75)(59.63,28.5)(55.5,27.25)
\qbezier(53,27)(48.75,26.13)(44.5,26.75)
\qbezier(41,27.75)(35.88,29.25)(32.25,31.75)
\qbezier(30.5,33)(28.25,34.75)(26,38.5)
\qbezier(24.75,40.5)(22.25,43)(21.75,46.5)
\qbezier(21.5,50)(20.75,51.25)(22,55.5)
\qbezier(22.5,58.25)(23.63,60.75)(24.25,62.25)
\put(16.5,47){\makebox(0,0)[cc]{$A_1$}}
\put(43.25,81){\makebox(0,0)[cc]{$B_2$}}
\put(68.5,76.5){\makebox(0,0)[cc]{$B_1$}}
\put(69.5,50.25){\makebox(0,0)[cc]{$C_2$}}
\put(22.25,45.25){\circle*{2.06}} \put(37.75,75){\circle*{2.5}}
\put(64.75,71.75){\circle*{1.58}} \put(73.75,51.25){\circle*{1.58}}
\put(60,28.75){\circle*{2.24}} \put(37,29){\circle*{2.12}}
\put(75.5,65.25){\makebox(0,0)[cc]{$\gamma H_i$}}
\put(35.25,25.25){\makebox(0,0)[cc]{$A_2$}}
\put(57.75,33.25){\makebox(0,0)[cc]{$C_1$}}
\end{picture}
\end{figure}

Consider a geodesic $\q_g$ in $\calg (p)$, and an arbitrary geodesic
triangle with $\q_g$ as an edge and the other two edges $\q_h\in
\calg(r_1)$ and $h\q_k$, where $\q_k\in \calg(r_2)$. Such a triangle
corresponds to a decomposition $g=hk$, where $h\in S_L(r_1)$ and
$k\in S_L (r_2)$. Let us apply property (*) to this geodesic
triangle. There exists $\gamma H_i\, ,\, i\in \{ 1,2,\dots m\}$,
such that $\overline{\nn}_\sigma (\gamma H_i)$ intersects the three
edges of the triangle. We denote as in the picture the respective
entrance and exit points of the edges into $\overline{\nn}_\sigma
(\gamma H_i)$ by $A_1, A_2, B_1, B_2, C_1, C_2$.

We can write $g=g_1\eta g_2$, where $g_1\in \gamma H_i$, $g_1$ is at
distance at most $\sigma$ from $A_2$, $\eta \in H_i$ and $g_1 \eta$
is at distance at most $\sigma $ from $C_1$. Property (*) implies
that $g_1$ is at distance at most $\delta+\sigma$ from $A_1$ and
that $g_1 \eta$ is at distance at most $\delta+\sigma$ from $C_2$.
It follows that there exists $\eta'\in H_i$ such that $g_1\eta'$ is
at distance at most $\sigma $ from $B_2$ and at distance at most
$\delta+\sigma $ from $B_1$. If we denote by $\eta''$ the element
$(\eta')\iv \eta$ then we have that $h=g_1\eta' g_3$ and that
$k=g_3\iv \eta''g_2$.



\me

\noindent \textit{Notation}: We denote by $\kappa$ the constant
$\sigma +\delta$.

\me

The previous considerations justify the following notations and
definitions.

\me

\noindent \textit{Notations}: We denote by $\Delta$ the set of all
$(g_1,g_2,g_3, \eta, \eta', \eta '')\in G^3 \times \bigsqcup_{i=1}^m
(H_i)^3$ such that:
    \begin{itemize}
        \item[(1)] $g=g_1 \eta g_2 \in S_L(p)$, $h= g_1 \eta' g_3\in
        S_L(r_1)$ and $k= g_3\iv \eta'' g_2 \in S_L (r_2)$;
        \item[(2)] $\eta =\eta' \eta ''$;
        \item[(3)] If $(\eta, \eta', \eta '')\in (H_i)^3$ for some $i\in
        \{1,2,\dots , m\}$, then the following hold:
\begin{itemize}
    \item $g_1$ is at distance at most $\kappa$ from the entrance point $A_2$
    of
$\q_g$ into $\overline{\nn}_\sigma (g_1
       H_i)$ and from the entrance point $A_1$ of $\q_h$ into
       $ \overline{\nn}_\sigma (g_1
       H_i)$;
    \item $g_1\eta$ is at distance at most $\kappa$ from the exit point $C_1$ of
    $\q_g$ from $ \overline{\nn}_\sigma (g_1
       H_i)$ and from the exit point $C_2$ of $\q_k$ from $\overline{\nn}_\sigma (g_1
       H_i)$;
    \item $g_1\eta'$ is at distance at most $\kappa$ from the exit point $B_2$ of
    $\q_h$ from $\overline{\nn}_\sigma (g_1
       H_i)$ and from the entrance point $B_1$ of $\q_k$ into $\overline{\nn}_\sigma (g_1
       H_i)$.
\end{itemize}
\end{itemize}
\me

\begin{definition}
 The set $\Delta$ is called the {\textit{set of central decompositions}}
  of geodesic triangles with edges in $\calg (p ) \times \calg (r_1) \times \calg
  (r_2)$.
  \end{definition}

  \me

\begin{definitions}
Given a geodesic $\q_g$ in $\calg(r)$ with $r \in \{ p \, ,\, r_2
\}$, we call every decomposition $g=g_1\eta g_2$ of $g$
corresponding to a central decomposition in $\Delta$ a
\textit{central decomposition of} $g$. We call $g_1$ and $g_2$ the
\textit{left} and \textit{right parts} of the decomposition; $\eta$
is called the \textit{middle part} of the decomposition of $g$.

For the fixed $g$ (and $\q_g$) we denote by $\call_g$ the set of
left parts of central decompositions of $g$, by $\calr_g$ the set
of right parts of central decompositions of $g$ and by $\ddd_g$
the set of triples $(g_1,\eta ,g_2)$ corresponding to central
decompositions of $g$. We also denote by $\call\calr_g$ the set of
pairs of left and right parts that can appear in a central
decomposition of $g$.
\end{definitions}

\me

We denote by $\ddd_p$ the set of all $\ddd_g$ with $g\in S_L(p)$.
The sets $\call\calr_p$, $\call_p$, $\calr_p$ are defined similarly.

\me

\noindent \textit{Notations}: Let $C$ be a subset of
$\prod_{i=1}^n X_i$. Let $a_i$ be a point in the $i$-th projection
of $C$, let $I$ be a subset in $\{ 1,2,\dots ,n\}\setminus \{ i\}$
and let $X^I=\prod_{i\in I} X_i$. We denote by $C^I(a_i)$ the set
of elements $\bar{x}$ in $X^I$ such that $(\bar{x},a_i)$ occurs in
the projection of $C$ in $X^I \times X_i$. Whenever there is no
risk of confusion, we drop the index $I$ in $C^I(a_i)$.

\me

For every fixed decomposition $d=(g_1,\eta ,g_2) \in \ddd_g$ we
also denote by $C_d$ the set $\Delta (d)$, that is, the set of
triples $(\eta',\eta'',g_3)$ such that $(g_1 ,g_2, g_3,\eta ,\eta'
,\eta '')$ is in $\Delta$. We denote by $\uuu_d$ the set of
projections on the last component of the triples in $C_d$, by
$\eee_d$ the set of projections $(\eta' ,\eta'' )$ on the first
two components of the triples in $C_d$ and by $E'_d$ and $E_d''$
the respective sets of $\eta'$ and $\eta''$.

\begin{lemma}\label{cdecomp}
Every $g\in S_L(r)$, with $r\in \{p\, ,\, r_2\}$, has at most
$C_1r_1+C_2$ central decompositions, where $C_1$ and $C_2$ are
universal constants.
\end{lemma}

\proof Suppose that $g_1$ is the left part of a central
decomposition of $g$. Then $g_1$ is at distance at most $\kappa $
from a point in a geodesic of length $r_1$. It follows that $g_1$
has length at most $r_1+\kappa $. On the other hand $g_1$ is at
distance at most $\kappa$ from a vertex $g_1'$ on $\q_g$, whence
$L(g_1')\leq r_1+2\kappa $. Thus, we have at most $r_1+2\kappa$
possibilities for $g_1'$. For each such $g_1'$ the number of left
cosets $\gamma H_i$ at distance at most $\sigma$ from it is bounded
by an universal constant, so $g_1'$ can be the entrance point of
$\q_g$ in $\overline{\nn}_\sigma (\gamma H_i)$ for at most a
constant number of left cosets $\gamma H_i$. The exit point $g_2'$
of $\q_g$ from $\overline{\nn}_\sigma (\gamma H_i)$ is uniquely
defined each time the left coset is fixed. For each left coset, the
number of points in it at distance at most $\sigma$ from $g_1'$ is
bounded by a universal constant, and likewise for $g_2'$. Thus, the
number of possibilities for $g_1$ and $g_1\eta$, once $g_1'$ is
fixed, is bounded by a constant. We deduce that on the whole there
are at most $K(r_1+2\kappa )$ possible central decompositions of $g$,
for some constant $K$.
\endproof

We continue the proof of Theorem \ref{main}. We can write that
$$
\| (x*y)_{p}\|^2 = \sum_{g\in S_L(p)} [(x*y)(g)]^2 =\sum_{g\in
S_L(p)} \left[ \sum_{(h,k)\in S_L(r_1)\times S_L(r_2),hk=g}x(h)y(k)
\right]^2 \leq $$
\begin{equation}\label{last}
\sum_{g\in S_L(p)} \left[ \sum_{d\in \ddd_g}
\sum_{(\eta',\eta'',g_3)\in C_d}x(g_1\eta' g_3)y(g_3\iv \eta'' g_2)
\right]^2\,.
\end{equation}

We use the following easy consequence of the Cauchy-Schwartz
inequality:
\begin{equation}\label{quadr}
\left( \sum_{i=1}^n a_i \right)^2 \leq n \sum_{i=1}^n a_i^2\, .
\end{equation}

This inequality, applied to the first sum in the brackets in the
last term of (\ref{last}), together with Lemma \ref{cdecomp}, gives
\begin{equation}\label{ineg1}
\| (x*y)_{p}\|^2 \leq (C_1r_1+C_2) \sum_{g\in S_L(p)} \sum_{d\in
\ddd_g} \left[ \sum_{(\eta',\eta'',g_3)\in C_d}x(g_1\eta'
g_3)y(g_3\iv \eta'' g_2) \right]^2\, .
\end{equation}

We re-write the term in the brackets on the right hand side of
inequality
 (\ref{ineg1}) and we apply the Cauchy-Schwartz
inequality as follows:

\begin{equation}\begin{array}{r} \sum_{(\eta',\eta'',g_3)\in C_d}x(g_1\eta'
g_3)y(g_3\iv \eta'' g_2) = \sum_{e=(\eta',\eta'')\in \eee_d }
\sum_{g_3\in C_d (e)} x(g_1\eta' g_3)y(g_3\iv \eta'' g_2)\leq\\
\label{xy} \sum_{e=(\eta',\eta'')\in \eee_d } \left[\sum_{g_3\in C_d
(e)} (x(g_1\eta' g_3))^2\right]^{1/2}\left[\sum_{g_3\in C_d (e)}
(y(g_3\iv \eta'' g_2))^2\right]^{1/2}\, .
\end{array}\end{equation}

We now define, for every $g\in S_L(p)$, $ \bar{g}=(g_1,g_2)\in \call
\calr_g$, $i\in \{1,2,\dots ,m \}$ and for $\eta =g_1\iv g g_2\iv
\in H_i\, $, the function $X_{\bar{g}}^i:H_i \to \R_+$ by
$$
X_{\bar{g}}^i (\eta')=\left[\sum_{g_3\in C_{({\bar{g}},\eta
)}(\eta' \, ,\, (\eta')\iv \eta )} (x(g_1\eta'
  g_3))^2\right]^{1/2}\mbox{ for every }\eta'\in
E'_{({\bar{g}},\eta )}\, ,
$$ and $X_{\bar{g}}^i (\eta')=0$ for all the other $\eta'\in H_i$.

Likewise we define the function $Y_{\bar{g}}^i:H_i \to \R_+$
by
$$
Y_{\bar{g}}^i (\eta'')=\left[\sum_{g_3\in C_{({\bar{g}},\eta )}(\eta
    (\eta'')\iv \, ,\, \eta'' )} (y(g_3\iv\eta''
  g_2))^2\right]^{1/2}\mbox{ for every }\eta''\in
E''_{({\bar{g}},\eta )}\, ,
$$ and $Y_{\bar{g}}^i (\eta'')=0$ for all the other $\eta''\in H_i$.

Then the sum in (\ref{xy}), if the middle part $\eta$ of $d$ is in
$H_i$, can be written as
$$
\sum_{(\eta' , \eta'' )\in  \eee_d \cap H_i \times H_i}
X_{\bar{g}}^i (\eta')Y_{\bar{g}}^i (\eta'')\, .
$$

Therefore inequality (\ref{ineg1}) gives

$$
\| (x*y)_{p}\|^2 \leq $$
\begin{equation}\label{2sums}(C_1r_1+C_2) \sum_{g\in
S_L(p)} \sum_{\bar{g}\in \call\calr_g}  \sum_{i=1}^m
\sum_{\eta \in \ddd_g (\bar{g})\cap H_i} \left[ \sum_{(\eta' ,
\eta'' )\in \eee_{(\bar{g},\eta )}\cap H_i \times H_i} X_{\bar{g}}^i
(\eta')Y_{\bar{g}}^i (\eta'') \right]^2 \, .
\end{equation}

We denote by $S$ the sum in the second term of the inequality
(\ref{2sums}), without the factor $C_1 r_1+C_2$. We have

$$\begin{array}{l}
S\leq \sum_{g\in S_L(p)} \sum_{\bar{g}\in \call\calr_g}
\sum_{i=1}^m \| X_{\bar{g}}^i * Y_{\bar{g}}^i \|^2 \leq \sum_{g\in
S_L(p)} \sum_{\bar{g}\in \call\calr_g} \sum_{i=1}^m
P_i(r_1+2\kappa)^2\| X_{\bar{g}}^i \|^2 \| Y_{\bar{g}}^i \|^2 \leq\\ \\
P(r_1) \sum_{g\in S_L(p)} \sum_{\bar{g}\in \call\calr_g}
\sum_{i=1}^m \| X_{\bar{g}}^i \|^2 \| Y_{\bar{g}}^i\|^2\, .
\end{array}
$$

The latter sum is smaller than
$$\begin{array}{l}
 \sum_{\bar{g}\in \call\calr_p} \left[\sum_{\eta \in
    \ddd_p (\bar{g}),\eta' \in E'_{(\bar{g},\eta )}} \sum_{g_3\in
    C_{(\bar{g},\eta )}(\eta', (\eta')\iv \eta )} (x(g_1\eta'
    g_3))^2\right]\cdot\\\hskip 2 in  \left[\sum_{\eta \in
    \ddd_p (\bar{g}),\eta'' \in E''_{(\bar{g},\eta )}} \sum_{g_3\in
    C_{(\bar{g},\eta )}(\eta(\eta'')\iv , \eta'' )} (y(g_3\iv\eta''
    g_2))^2\right]\leq \\

 \left[\sum_{g_1\in \call_p} \sum_{(\eta', g_3
)\in
    \Delta (g_1)} (x(g_1\eta'
    g_3))^2\right]\cdot \left[\sum_{g_2\in \calr_p } \sum_{(g_3,\eta'')\in
    \Delta (g_2)} (y(g_3\iv\eta''
    g_2))^2\right]\, .
\end{array}$$

Every $x(g)^2$ with $g\in S_L(r_1)$ appears at most $K_1r_1$ times
in the first sum above, where $K_1$ is an universal constant, hence
the first sum is at most $K_1r_1 \| x\|^2$. Lemma \ref{cdecomp}
applied to every $g\in S_L (r_2)$ implies that every $y(g)^2$ with
$g\in S_L(r_2)$ appears at most $C_1r_1+C_2$ times in the second sum
above, hence the second sum is at most $(C_1r_1+C_2) \| y\|^2$.

We deduce that
$$
S\leq Q_1(r_1)P(r_1)\| x\|^2 \| y\|^2\, ,
$$ where $Q_1$ is an universal polynomial of degree $2$. We may conclude that
$$
\| (x*y)_p \|^2 \leq Q(r_1)P(r_1)\| x\|^2 \| y\|^2\, ,
$$ where $Q$ is an universal polynomial of degree $3$.
\endproof

\begin{remark} \label{twostars} Theorem \ref{main} holds if we
replace property (*) in the definition of (*)-relative hyperbolicity
by the following weaker property.

We say that a countable group $G$ endowed with a length-function $L$ is
(**)-relatively hyperbolic with respect to the subgroups $H_1,...,H_m$
if there exist two polynomials $Q_1(r)$ and $Q_2(r)$ such that the
following condition holds.

\begin{quote}(**) Let ${\mathcal{H}}=\bigcup_{i=1}^m (H_i\times H_i)$. There exists a map $T\colon G\times
G\to G\times {\mathcal{H}}$ such that

(i) if $T(g,h)=(a,g',h')$ then $T(h,g)=(a,h',g')$ and $$T(h\iv, h\iv
g)=(h\iv ah', (h')\iv, (h')\iv g');$$

(ii) $L(h')\le Q_1(L(h))$;

(iii) for each $g\in G$ the number of pairs $(a,g')$ such that $(a,g',\, \cdot \, )=T(g,h)$ for
some $h$ with $L(h)=r$ does not exceed $Q_2(r)$.
\end{quote}

One can interpret a pair $(g,h)$ as two sides of the triangle with
vertices $x, xg, xh$ (for any basepoint $x\in G$). The image
$T(g,h)$ gives three vertices of a triangle $xa, xag', xah'$ from a
coset of $H_i$. Condition (i) means that the set $\{a, ag', ah'\}$
is stable under permutations of the vertices $x, xg, xh$. Condition
(ii) means that the sides of the resulting triangle are not too
large. The pair $(a,g')$ can be interpreted as a central
decomposition of $g$ corresponding to the triangle $(1, g, h)$ (more
precisely the left and the middle part of the decomposition, the
right part is uniquely determined by the other two). Condition (iii)
means that the number of central decompositions of $g$ corresponding
to triangles $(x, xg, xh)$ with $L(h)=r$ is bounded by a polynomial
in $r$. Condition (iii) replaces Lemma \ref{cdecomp}, the proof of
Theorem \ref{main} carries almost verbatim.

The class of (**)-relatively hyperbolic groups is certainly wider
than the class of relatively hyperbolic groups. For example $\Z^2$
is (**)-relatively hyperbolic with respect to the trivial subgroup.
To construct a map $T$, consider the two $\Gamma$-shaped geodesics in
$\Z^2$ connecting any two points $A,B$. For every pair
$\{g,h\}$ choose one of the two possible geodesics for each pair
$(1,g)$, $(1,h)$, $(g,h)$ so that the three geodesics intersect in a
point $a$. The map $T$ takes $(g,h)$ to $(a,1,1)$. One should choose
the point $a$ so that condition (i) is satisfied: if we choose a certain
$a$ for the pair $(g,h)$ then the corresponding $a'$s for $(h,g)$ and
$(h\iv, h\iv g)$ are determined uniquely. Then the conditions (ii) and
(iii) are also satisfied, for $Q_1(r)=0$, $Q_2(r)=2r$.

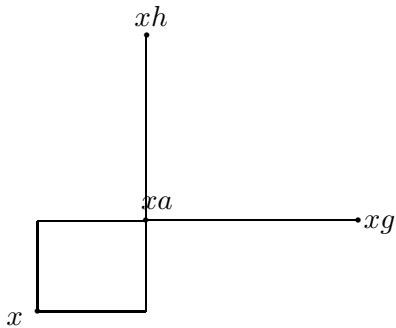
\begin{center}
\unitlength .5mm 
\linethickness{0.4pt}
\ifx\plotpoint\undefined\newsavebox{\plotpoint}\fi 
\begin{picture}(115.5,98.75)(0,0)


\put(26.5,14){\line(1,0){28.75}}

\put(55.25,14){\line(0,1){24}}

\put(26.5,14){\line(0,1){24}}

\put(26.5,38){\line(1,0){28.75}}

\multiput(55.25,14)(0,9.1875){8}{\line(0,1){9.1875}}
\put(55.25,38.25){\line(1,0){56.75}}
\put(26.5,14){\circle*{1.58}} \put(55.25,38.25){\circle*{1.58}}
\put(55.5,87.5){\circle*{1.58}} \put(111.75,38.25){\circle*{1.58}}
\put(20.5,12){\makebox(0,0)[cc]{$x$}}
\put(56.75,92.75){\makebox(0,0)[cc]{$xh$}}
\put(117.5,37){\makebox(0,0)[cc]{$xg$}}
\put(58.25,43){\makebox(0,0)[cc]{$xa$}}
\end{picture}
\end{center}
Similarly, for every $k$, $\Z^k$ is (**)-relatively hyperbolic with
respect to the trivial subgroup.

More generally, every group with polynomial growth function $f(n)$ is
(**)-relatively hyperbolic with respect to the trivial subgroup: for
every pair $g,h$ pick a shortest side of the
triangle with vertices $1, g, h$. Let $a$ be $g$ if $g$ is a
vertex of the chosen shortest side, or $1$ otherwise (again choose $a$
depending on $(g,h)$ so that (i) is satisfied). The map $T$ takes
$(g,h)$ to $(a,1,1)$. It is an easy exercise to check that
conditions (ii) and (iii) are satisfied with $Q_1(r)=0$,
$Q_2(r)=2f(r)+2$.

A group $G$ that is (*)-relatively hyperbolic (in particular
relatively hyperbolic) with respect to the subgroups $H_1,...,H_m$ is
(**)-relatively hyperbolic with respect to these subgroups, by Lemma
\ref{cdecomp} (both $Q_1$ and $Q_2$ are linear polynomials).

It would be interesting to see how large the class of
(**)-relatively hyperbolic groups is. In particular, it would be
interesting to know if every mapping class group is (**)-relatively
hyperbolic with respect to direct products of mapping class groups
of smaller genus (in that case we would be able to proceed by
induction and prove that all mapping class groups have RD).
\end{remark}

\addtocontents{toc}{\contentsline {section}{\numberline { }
References \hbox {}}{\pageref{bibbb}}}

\bigskip

\begin{minipage}[t]{2.9 in}
\noindent Cornelia Dru\c tu\\ Department of Mathematics\\
University of Lille 1\\ and UMR CNRS 8524 \\ Cornelia.Drutu@math.univ-lille1.fr\\
\end{minipage}
\begin{minipage}[t]{2.6 in}
\noindent Mark V. Sapir\\ Department of Mathematics\\
Vanderbilt University\\
msapir@math.vanderbilt.edu\\
\end{minipage}

\end{document}